\documentclass{gtart}

\def\ifplaintex{\expandafter\ifx\csname documentclass\endcsname\relax}

\ifplaintex 
\else
\fi

\def\gtm{{\mathsurround=0pt\it $\cal G\mskip-2mu$eometry \&\ 
$\cal T\!\!$opology $\cal M\mskip-1mu$onographs}}    %  for monographs

\def\gtp{{\mathsurround=0pt\it $\cal G\mskip-2mu$eometry \&\ 
$\cal T\!\!$opology $\cal P\!$ublications}}  % GT publications

\def\recd{{\small Received:\qua\receiveddate\ifx\reviseddate\relax
\else\qquad Revised:\qua\reviseddate\fi\par}} 

\def\volumenumber#1{\def\thevolumenumber{#1}}
\def\volumeyear#1{\def\thevolumeyear{#1}}
\def\volumename#1{\def\thevolumename{#1}}
\def\papernumber#1{\def\thepapernumber{#1}}
\def\pagenumbers#1#2{\def\startpage{#1}\def\finishpage{#2}}
\def\published#1{\def\publishdate{#1}}
\def\received#1{\def\receiveddate{#1}}

\def\accepted#1{\def\accepteddate{#1}}
\def\asciititle#1{\def\theasciititle{#1}}

\def\coverauthors#1{\def\thecoverauthors{#1}}
\def\asciiauthors#1{\def\theasciiauthors{#1}}
\def\asciiaddress#1{\def\theasciiaddress{#1}}

\def\coverauthors#1{\def\thecoverauthors{#1}}
\long\def\asciiabstract#1{\long\def\theasciiabstract{#1}}

\def\shorttitle#1{\def\theshorttitle{#1}}

\let\\\par
\let\thevolumenumber\relax\let\thepapernumber\relax
\let\thevolumeyear\relax\let\startpage\relax
\let\finishpage\relax\let\publishdate\relax\let\receiveddate\relax
\let\reviseddate\relax\let\accepteddate\relax\let\theasciititle\relax
\let\theasciiauthors\relax\let\theasciiaddress\relax
\let\theasciiabstract\relax
\let\thecoverauthors\relax
\let\thecoverauthors\relax\let\theerratum\relax\let\theasciiemail\relax
\let\theshortauthors\relax\let\theshorttitle\relax

\def\startpage{1}\def\finishpage{15}\def\thepapernumber{77}

\volumenumber{2}
\volumename{Proceedings of the Kirbyfest}
\volumeyear{1999}

\long\def\maketitlep{   % start of definition of \maketitlep

\count0=\startpage

\gtm\nl        %   GT mongraphs (top left) 
{\small Volume \thevolumenumber: \thevolumename\nl 
\ifx\theerratum\relax\else Erratum \erratumnumber\nl\fi
Pages \startpage--\finishpage\nl}

\vglue 0.1truein   % top margin

{\parskip=0pt\leftskip 0pt plus 1fil\def\\{\par\smallskip}{\ifplaintex\large
\else\Large\fi\bf\thetitle}\par\medskip}   
\vglue 0.05truein 

{\parskip=0pt\leftskip 0pt plus 1fil\def\\{\par}{\sc\theauthors}
\par\medskip}%
 
\vglue 0.03truein 

{\small\leftskip 25pt\rightskip 25pt{\bf Abstract}\stdspace\theabstract

{\bf AMS Classification}\stdspace\theprimaryclass
\ifx\thesecondaryclass\relax\else; \thesecondaryclass\fi\par
{\bf Keywords}\stdspace \thekeywords\par}\vglue 7pt

}   % end of definition of \maketitlep

\font\phead=cmsl9 scaled 950
\font=cmsl9 scaled 1050
\font\pnum=cmbx10 scaled 913
\font\lnum=cmbx10 
\font\pfoot=cmsl9 scaled 950
\font=cmsl9 scaled 1050
\ifplaintex
\headline{\vbox to 0pt{\vskip -4.5mm\line{\small\phead\ifnum
\count0=\startpage ISSN 1464-8997 (on line)
1464-8989 (printed) \hfill {\pnum\folio}\else\ifodd\count0\def\\{ }% 
\ifx\theshorttitle\relax\thetitle\else\theshorttitle\fi\hfill{\pnum\folio}
\else\def\\{ and }{\pnum\folio}\hfill\ifx\theshortauthors\relax\theauthors
\else\theshortauthors\fi\fi\fi}\vss}}
\footline{\vbox to 0pt{\vglue 0mm\line{\small\pfoot\ifnum\count0=\startpage
Published \publishdate:\qua\copyright\ \gtp\hfill\else
\gtm, Volume \thevolumenumber\ (\thevolumeyear)\hfill\fi}\vss
}}
\else
\makeatletter
\def\@oddhead{{\small\lhead\ifnum\count0=\startpage ISSN 1464-8997 (on line)
1464-8989 (printed) \hfill {\lnum\number\count0}\else\ifodd\count0
\def\\{ }\ifx\theshorttitle\relax \thetitle \else\theshorttitle\fi\hfill
{\lnum\number\count0}\else\def\\{ and }{\lnum\number\count0}
\hfill\ifx\theshortauthors\relax 
\theauthors\else\theshortauthors\fi\fi\fi}}\def\@evenhead{@oddhead}
\def\@oddfoot{\small\lfoot\ifnum\count0=\startpage Published \publishdate:\qua\copyright\ \gtp\hfill\else
\gtm, Volume \thevolumenumber\ (\thevolumeyear)\hfill\fi}
\def\@evenfoot{@oddfoot}
\makeatother
\fi

\let\maketitlepage\maketitlep
\let\makeshorttitle\maketitlepage
\let\maketitle\maketitlepage

\newwrite\gtoutfile
\long\gdef\makeheadfile{  %%% start of definition of \makeheadfile
{\def\\{, }\def\s{ }
\immediate\openout\gtoutfile head.xxx
\immediate\write\gtoutfile{To: math@arxiv.org}
\immediate\write\gtoutfile{Subject: put OR rep NNNNN:ppppp}
\immediate\write\gtoutfile{--text follows this line--}
\immediate\write\gtoutfile{Proxy-for: \ifx\theasciiauthors\relax
\theauthors\else\theasciiauthors\fi\s<\ifx\theasciiemail\relax\theemail\else\theasciiemail\fi>}
\immediate\write\gtoutfile{\noexpand\\}
\immediate\write\gtoutfile{Authors: \ifx\theasciiauthors\relax
\theauthors\else\theasciiauthors\fi}
{\def\\{ }\immediate\write\gtoutfile{Title: \ifx\theasciititle\relax
\thetitle\else\theasciititle\fi}}
\immediate\write\gtoutfile{Subj-class: GT or SG, GR etc}
\immediate\write\gtoutfile{MSC-class: \theprimaryclass\ifx\thesecondaryclass\relax\else, \thesecondaryclass\fi}
\immediate\write\gtoutfile{Journal-ref: Geom. Topol. Monogr. \thevolumenumber\s
(\thevolumeyear) \startpage-\finishpage}
\immediate\write\gtoutfile{Comments: Published by Geometry and Topology Monographs at}
\immediate\write\gtoutfile{\s\s\s  http://www.maths.warwick.ac.uk/gt/GTMon\thevolumenumber/paper\thepapernumber.abs.html}
\immediate\write\gtoutfile{\noexpand\\}
\immediate\write\gtoutfile{}
\ifx\theasciiabstract\relax
\immediate\write\gtoutfile{\theabstract}\else
\immediate\write\gtoutfile{\theasciiabstract}\fi
\immediate\write\gtoutfile{}
\immediate\write\gtoutfile{\noexpand\\}
\immediate\write\gtoutfile{}
\immediate\closeout\gtoutfile}}  %%% end of definition of \makeheadfile

\def\maketitlepage{\maketitlep\makeheadfile}
\let\makeshorttitle\maketitlepage
\let\maketitle\maketitlepage

\volumenumber{4}
\volumename{Invariants of knots and 3-manifolds (Kyoto 2001)}
\volumeyear{2002}
\papernumber{6}
\pagenumbers{69}{87}
\received{3 December 2001}
\accepted{22 July 2002}
\published{19 September 2002}

\usepackage{amsmath,amssymb,texdraw}

\newcommand{\Z}{\mathbb{Z}}
\newcommand{\Q}{\mathbb{Q}}
\newcommand{\R}{\mathbb{R}}
\newcommand{\C}{\mathbb{C}}
\newcommand{\Io}{\mathbb I}

\newcommand{\ep}{\epsilon}
\newcommand{\vep}{\varepsilon}
\newcommand{\I}{\sqrt{-1}}
\newcommand{\npr}{|\triangle_+|}
\newcommand{\wl}{\Lambda^W}
\newcommand{\rl}{\Lambda^R}
\newcommand{\la}{\langle}
\newcommand{\ra}{\rangle}
\newcommand{\ria}{\rightarrow}

\newcommand{\vol}{\operatorname{vol}}
\newcommand{\sign}{\operatorname{sign}}
\newcommand{\inte}{\operatorname{int}}

\newcommand{\id}{\text{\rm id}}
\newcommand{\Hom}{\text{\rm Hom}}
\newcommand{\ns}{\text{\rm n}}
\newcommand{\os}{\text{\rm o}}
\newcommand{\s}{\text{\rm s}}
\newcommand{\tO}{\text{\rm O}}
\newcommand{\tr}{\text{\rm tr}}

\newcommand{\frg}{\mathfrak g}
\newcommand{\frh}{\mathfrak h}
\newcommand{\frhR}{{\mathfrak h}_{\mathbb{R}}}
\newcommand{\frsl}{\mathfrak sl}

\newcommand{\mA}{\mathcal A}
\newcommand{\mC}{\mathcal C}
\newcommand{\mD}{\mathcal D}
\newcommand{\mR}{\mathcal R}
\newcommand{\mT}{\mathcal T}
\newcommand{\mV}{\mathcal V}

\newtheorem{thm}{Theorem}[section]
\newtheorem{conj}[thm]{Conjecture}
\newtheorem{cor}[thm]{Corollary}
\newtheorem{lem}[thm]{Lemma}
\newtheorem{prop}[thm]{Proposition}

\theoremstyle{definition}
\newtheorem{rem}[thm]{Remark}

\newcommand{\refthm}[1]{Theorem~\ref{#1}}
    
    \newcommand{\refcor}[1]{Corollary~\ref{#1}}
    \newcommand{\reflem}[1]{Lemma~\ref{#1}}
    \newcommand{\refprop}[1]{Proposition~\ref{#1}}
    \newcommand{\refrem}[1]{Remark~\ref{#1}}

\let\HS\qed

\begin{document}

\title{Quantum invariants of Seifert $3$--manifolds\\and their
asymptotic expansions} 

\asciititle{Quantum invariants of Seifert 3-manifolds and their
asymptotic expansions} 

\shorttitle{Quantum invariants of Seifert $3$--manifolds}

\author{S\o ren Kold Hansen\\Toshie Takata}
\coverauthors{S\noexpand\o ren Kold Hansen\\Toshie Takata}
\asciiauthors{Soren Kold Hansen\\Toshie Takata}

\address{School of Mathematics, University of Edinburgh, 
JCMB\\King's Buildings, Edinburgh EH9 3JZ, UK}
\email{hansen@maths.ed.ac.uk, takata@math.sc.niigata-u.ac.jp}
\secondaddress{Department of Mathematics, Faculty of Science\\Niigata University, Niigata 950-2181, Japan}

\asciiaddress{Dept of Maths and Stats, University of Edinburgh,
JCMB\\King's Buildings, Edinburgh EH9 3JZ, UK\\Department of
Mathematics, Faculty of Science\\Niigata University, Niigata 950-2181,
Japan}

\begin{abstract}
We report on recent results of the authors concerning calculations of
quantum invariants of Seifert $3$--manifolds. These results include a
derivation of the Reshetikhin--Turaev invariants of all oriented
Seifert manifolds associated with an arbitrary complex finite
dimensional simple Lie algebra, and a determination of the asymptotic
expansions of these invariants for lens spaces. Our results are in
agreement with the asymptotic expansion conjecture due to J.\ E.\
Andersen \cite{Andersen1}, \cite{Andersen2}.
\end{abstract}

\asciiabstract{We report on recent results of the authors concerning
calculations of quantum invariants of Seifert 3-manifolds. These
results include a derivation of the Reshetikhin-Turaev invariants of
all oriented Seifert manifolds associated with an arbitrary complex
finite dimensional simple Lie algebra, and a determination of the
asymptotic expansions of these invariants for lens spaces. Our results
are in agreement with the asymptotic expansion conjecture due to JE
Andersen [The Witten invariant of finite order mapping tori I, to
appear in J. Reine Angew. Math.] and [The asymptotic expansion
conjecture, from `Problems on invariants of knots and $3$--manifolds',
edited by T. Ohtsuki, http://www.ms.u-tokyo.ac.jp/\char'176tomotada/proj01/].}

\primaryclass{57M27}

\secondaryclass{17B37, 18D10, 41A60}

\keywords{Quantum invariants, Seifert manifolds, modular categories, quantum groups, asymptotic expansions}

%\maketitlepage
\makeshorttitle

\def\undersmile#1{\lower5.7pt\hbox{$\smallsmile$}\kern-0.6em #1}

\section{Introduction}\label{sec-Introduction}

In 1988 E.\ Witten \cite{Witten} proposed new invariants
$Z_{k}^{G}(X,L) \in \C$ of an
arbitrary closed oriented $3$--manifold $X$ with
an embedded colored link $L$ by
quantizing the Chern--Simons field theory associated to a
simple and simply connected compact Lie group $G$,
$k$ being an arbitrary positive integer called the (quantum) level.
The invariant $Z_{k}^{G}(X,L)$ is given by a Feynman path integral 
over the (infinite dimensional) space of gauge equivalence
classes of connections in a $G$ bundle over $X$.
This integral should be understood in a formal way since,
at the moment of
writing, it seems that no mathematically rigorous definition is known,
cf.\ \cite[Sect.~20.2.A]{JohnsonLapidus}.
The invariants $Z_{k}^{G}$ are called the quantum
$G$--invariants or Witten's invariants associated to
$G$.

Shortly afterwards, N.\ Reshetikhin and V.\ G.\ Turaev
\cite{ReshetikhinTuraev} defined in a mathematically rigorous
way invariants $\tau_{r}^{\frsl_{2}(\C)}(X,L) \in \C$
of the pair $(X,L)$
by combinatorial means using irreducible
representations of the quantum deformations of $\frsl_{2}(\C)$
at certain roots of unity, $r$ being an integer $\geq 2$ associated
to the order of the root of unity.
Later quantum invariants $\tau_{r}^{\frg}(X,L) \in \C$ associated to
other complex simple Lie algebras $\frg$ were constructed using 
representations of the quantum deformations of $\frg$ at
`nice' roots of unity, see \cite{TuraevWenzl}. 
We call $\tau_{r}^{\frg}$ for the quantum $\frg$--invariants or
the RT--invariants associated to $\frg$.

Both in Witten's approach and in the
approach of Reshetikhin and Turaev the invariants are part of a 
topological quantum field theory (TQFT) (or more correctly
a family of TQFT's). This implies that the
invariants are defined for compact oriented $3$--dimensional cobordisms 
(perhaps with some extra structure on the boundary),
and satisfy certain cut-and-paste axioms, see \cite{Atiyah}, \cite{Blanchetetal},
\cite{Quinn}, \cite{Turaev}. The TQFT's of Reshetikhin and Turaev can from
an algebraic point of view be given a more general formulation by using
so-called modular (tensor) categories \cite{Turaev}. The representation theory
of the quantum deformations of $\frg$ at certain roots of unity,
$\frg$ an arbitrary finite dimensional
complex simple Lie algebra, induces such modular categories,
see e.g.\ \cite{Kirillov}, \cite{BakalovKirillov}, \cite{Le}.

For an invariant to be powerful one should be able to calculate it.
A problem with the quantum invariants of knots and $3$--manifolds
is that they are rather hard to calculate. In fact people have only been
able to calculate these invariants for certain (families of)
knots and $3$--manifolds.
The lens spaces and more generally the Seifert $3$--manifolds constitute such a
family, and there is a wealth of literature about different
calculations of quantum invariants
of these spaces, see \cite[Introduction]{Hansen1} for some references.
In \cite{Hansen1} the RT--invariant associated to an
arbitrary modular category is calculated for any Seifert manifold,
cf.\ \cite[Theorem 4.1]{Hansen1}.
(Here and in the rest of this paper a $3$--manifold means a closed oriented
$3$--manifold. In particular, a Seifert manifold is an oriented Seifert manifold.)

A solution to the above problem and to the general
problem of understanding the
topological `meaning' of
the quantum invariants could be to
determine relationships between the
quantum invariants and classical (well understood and calculable) invariants.
However, this seems to be a rather hard task.
This leads us into one of the themes in this article,
namely asymptotic expansions of the invariants. 
By using stationary phase approximation techniques together with path 
integral arguments Witten was able \cite{Witten}
to express the leading asymptotics 
of $Z_{k}^{G}(X)$ in the limit $k \ria \infty$
as a sum over the set
of stationary points for the Chern--Simons functional.
The terms in this sum are expressed by such 
topological/geometric invariants as Chern--Simons invariants, Reidemeister
torsions and spectral flows, so here we see a way to extract topological
information from the invariants. 
A full asymptotic expansion of Witten's invariant
is expected on the basis of a full perturbative 
analysis of the Feynman path integral, see \cite{AxelrodSinger1},
\cite{AxelrodSinger2}, \cite{Axelrod}.
It is generally believed that the family of TQFT's of Reshetikhin and Turaev is
a mathematical realization of Witten's family of TQFT's. This belief has together
with known results concerning asymptotics of the RT--invariants lead
to a conjecture, the asymptotic expansion conjecture (AEC),
which specifies in a rather precise way the asymptotic behaviour of the
RT--invariants.
The AEC was proposed by Andersen in \cite{Andersen1},
where he proved it for mapping 
tori of finite order diffeomorphisms of orientable surfaces of
genus at least two
using the gauge theoretic definition of the quantum invariants.

In this paper we explain recent results of the authors concerning the
RT--invariants of Seifert manifolds. Explicitly we state formulas
for the invariants $\tau_{r}^{\frg}$ of all Seifert manifolds in terms
of the Seifert invariants and standard data for $\frg$, $\frg$ being
an arbitrary complex finite dimensional simple Lie algebra,
cf.\ \refthm{Lie-Seifert}. Moreover,
we analyse more carefully the invariants $\tau_{r}^{\frg}(X)$ for $X$
any lens space, thereby determining a formula for
the large $r$ asymptotics of these invariants, cf.\ \refthm{asymp-lens}
and the remark following this theorem.
This formula is in agreement with the AEC.

A part of the paper is concerned with studying
a certain family of 
finite dimensional complex representations $\mR=\mR_{r}^{\frg}$ 
of $SL(2,\Z)$. These representations
are known from the study of theta functions and modular forms in
connection with the study of affine Lie algebras, cf.\ \cite{KacPeterson},
\cite[Sect.~13]{Kac}.
They also play a fundamental role in
conformal field theory and (therefore) in the Chern--Simons TQFT's of Witten,
see e.g.\ \cite{GepnerWitten}, \cite{Verlinde}, \cite{Witten}.
In case $\frg=\frsl_{2}(\C)$, Jeffrey \cite{Jeffrey1}, \cite{Jeffrey2} has
determined a nice formula for $\mR_{r}^{\frg}(U)$ in terms of the entries in
$U \in SL(2,\Z)$. \refthm{rep} is a direct
extension of Jeffrey's result to arbitrary $\frg$.
The representations $\mR_{r}^{\frg}$ are of interest when calculating 
the RT--invariants of the Seifert manifolds since certain matrices,
which can be expressed through $\mR_{r}^{\frg}$, enter into the formulas of the invariants.

The paper is organized as follows. In Sect.~\ref{sec-Seifert-manifolds}
we introduce notation for the Seifert manifolds and recall surgery
presentations for these manifolds. In 
Sect.~\ref{sec-Quantum-invariants}
we explain our calculation of the $\frg$--invariants of the Seifert manifolds.
In Sect.~\ref{sec-The-asymptotic}
we state the asymptotic
expansion conjecture and determine the asymptotic expansions
of the $\frg$--invariants of the lens spaces.
In the appendix we sketch the proof of the formula for the entries
in $\mR_{r}^{\frg}(U)$, $U \in SL(2,\Z)$, \refthm{rep}.
The paper is to some extend expository.
Details and most technicalities, in particular in connection to
the proof of \refthm{rep}, will be given in \cite{HansenTakata}.

\rk{Acknowledgements} This work were done while the first author was
supported by a Marie Curie Fellowship of the European Commission (CEE
$\text{N}^{o}$ HPMF--CT--1999--00231).
He acknowledge hospitality of
l'Institut de Recherche Math\'{e}matique Avanc\'{e}e (IRMA), Universit\'{e}
Louis Pasteur and C.N.R.S., Strasbourg, while being a Marie Curie Fellow.
A part of this work was done while the second author visited IRMA.
She thanks this department for hospitality during here stay.
Another part was done while both authors visited the Research Institute
for Mathematical Sciences (RIMS), Kyoto University.
We would like to thank RIMS for hospitality during the special month on 
{\it Invariants of knots and $3$--manifolds}, September 2001.
We also thank the organisers of this program for letting
us present this work at the workshop of the special month.
Finally the first author thanks J.\ E.\ Andersen for helpful
conversations about quantum
invariants in general and about asymptotics of these invariants in particular.

\section{Seifert manifolds}\label{sec-Seifert-manifolds}

For Seifert manifolds we will use the notation introduced by Seifert in his
classification results for these manifolds, see \cite{Seifert1},
\cite{Seifert2}, \cite[Sect.~2]{Hansen1}. That is, 
$(\ep;g \, | \, b;(\alpha_{1},\beta_{1}),\ldots,(\alpha_{n},\beta_{n}))$ is
the Seifert manifold with orientable base of genus $g \geq 0$ if
$\ep=\os$ and non-orientable base of genus $g>0$ if $\ep=\ns$ 
(where the genus of the non-orientable connected sum $\# ^{k} \R \text{P}^{2}$
is $k$). (In \cite{Seifert1}, \cite{Seifert2}
$(\ep;g \, | \, b;(\alpha_{1},\beta_{1}),\ldots,(\alpha_{n},\beta_{n}))$ is
denoted $(\tO,\ep;g \, | \, b;\alpha_{1},\beta_{1};\ldots;\alpha_{n},\beta_{n})$,
but we leave out the $\tO$, since we are only dealing with
oriented Seifert manifolds.)
The pair $(\alpha_{j},\beta_{j})$ of coprime integers is the (oriented)
Seifert invariant of the $j$'th exceptional (or singular) fiber.
We have $0< \beta_{j}<\alpha_{j}$. The integer $-b$ is equal to the
Euler number of the Seifert fibration $(\ep;g \, | \, b)$
(which is a locally trivial $S^{1}$--bundle).
More generally, the Seifert Euler number of
$(\ep;g \, | \, b;(\alpha_{1},\beta_{1}),\ldots,(\alpha_{n},\beta_{n}))$
is $E=-\left(b+\sum_{j=1}^{n} \beta_{j}/\alpha_{j} \right)$.
We note that lens spaces are Seifert manifolds with base $S^{2}$
and zero, one or two exceptional fibers.
According to \cite[Fig.~12 p.~146]{Montesinos}, the manifold
$( \ep ;g \, | \, b;(\alpha_{1},\beta_{1}),\ldots,(\alpha_{n},\beta_{n}))$ 
has a surgery
presentation as shown in Fig.~\ref{fig-A1} if $\ep=\os$ and as
shown in Fig.~\ref{fig-A2} if $\ep=\ns$.
The $\undersmile{g}$ indicate $g$ repetitions.

\begin{figure}[ht!]

\begin{center}
\begin{texdraw}
\drawdim{cm}

\setunitscale 0.6

\linewd 0.02 \setgray 0

\move(0 4)

\move(0 0) \lellip rx:3 ry:1.6

\linewd 0.2 \setgray 1

\move(4 0) \larc r:1.5 sd:218 ed:200
\move(4 0) \larc r:3 sd:213 ed:204
\move(-4 0) \larc r:1.5 sd:-20 ed:150
\move(-4 0) \larc r:2 sd:-23 ed:150
\move(-4 0) \larc r:3.5 sd:-20 ed:155
\move(-4 0) \larc r:3 sd:-23 ed:155

\linewd 0.02 \setgray 0

\move(-4 0) \larc r:1.5 sd:-20 ed:150
\move(-4 0) \larc r:2 sd:-23 ed:150
\move(-4 0) \larc r:3.5 sd:-20 ed:155
\move(-4 0) \larc r:3 sd:-23 ed:155

\move(-5.88 0.67) \clvec(-5.85 0.9)(-5.6 0.95)(-5.45 0.8)
\move(-5.88 -0.67) \clvec(-5.85 -0.9)(-5.6 -0.95)(-5.45 -0.8)

\move(-5.3 0.75) \clvec(-5.35 0.55)(-5.6 0.5)(-5.75 0.7)
\move(-5.3 -0.75) \clvec(-5.35 -0.55)(-5.6 -0.5)(-5.75 -0.7)

\move(-7.3 1.17) \clvec(-7.25 1.4)(-7 1.45)(-6.85 1.3)
\move(-7.3 -1.17) \clvec(-7.25 -1.4)(-7 -1.45)(-6.85 -1.3)

\move(-6.7 1.29) \clvec(-6.75 1.09)(-7 1.04)(-7.15 1.24)
\move(-6.7 -1.29) \clvec(-6.75 -1.09)(-7 -1.04)(-7.15 -1.24)

\move(-4 0) \larc r:1.5 sd:210 ed:320
\move(-4 0) \larc r:2 sd:210 ed:318
\move(-4 0) \larc r:3.5 sd:205 ed:330
\move(-4 0) \larc r:3 sd:205 ed:328
\move(-4 0) \larc r:1.5 sd:160 ed:200
\move(-4 0) \larc r:2 sd:160 ed:200
\move(-4 0) \larc r:3.5 sd:160 ed:200
\move(-4 0) \larc r:3 sd:160 ed:200

\move(4 0) \larc r:1.5 sd:218 ed:200
\move(4 0) \larc r:3 sd:213 ed:204

\move(-6.8 -0.1) \htext{$\cdots$}
\move(1.5 -0.1) \htext{$\cdots$}
\move(-6.85 -0.8) \htext{$\undersmile{g}$}
\move(-6.65 2.8) \htext{$0$}
\move(-5.75 1.6) \htext{$0$}
\move(-7.95 0.3) \htext{$0$}
\move(-6.45 0.15) \htext{$0$}
\move(-3.95 -0.3) \htext{$-b$}
\move(5.1 1) \htext{$\frac{\alpha_{1}}{\beta_{1}}$}
\move(6.2 2.2) \htext{$\frac{\alpha_{n}}{\beta_{n}}$}

\move(0 -4)

\end{texdraw}
\end{center}

\caption{Surgery presentation of $(\os ;g \, | \, b;(\alpha_{1},\beta_{1}),\ldots,(\alpha_{n},\beta_{n}))$}\label{fig-A1}
\end{figure}

\begin{figure}[ht!]

\begin{center}
\begin{texdraw}
\drawdim{cm}

\setunitscale 0.6

\linewd 0.02 \setgray 0

\move(0 4)

\move(0 0) \lellip rx:3 ry:1.6
\move(-7.5 0) \lellip rx:0.5 ry:0.3
\move(-5.5 0) \lellip rx:0.5 ry:0.3

\linewd 0.2 \setgray 1
\move(-4 0) \larc r:1.5 sd:-20 ed:180
\move(-4 0) \larc r:3.5 sd:-20 ed:180
\move(4 0) \larc r:1.5 sd:215 ed:200
\move(4 0) \larc r:3 sd:212 ed:205

\linewd 0.02 \setgray 0

\move(-4 0) \larc r:1.5 sd:-20 ed:180
\move(-4 0) \larc r:3.5 sd:-20 ed:180

\move(-4 0) \larc r:1.5 sd:200 ed:320
\move(-4 0) \larc r:3.5 sd:188 ed:330

\move(4 0) \larc r:1.5 sd:215 ed:200
\move(4 0) \larc r:3 sd:212 ed:205

\move(-6.8 -0.1) \htext{$\cdots$}
\move(1.5 -0.1) \htext{$\cdots$}
\move(-6.8 -0.9) \htext{$\undersmile{g}$}
\move(-8.3 0.25) \htext{$2$}
\move(-6.3 0.25) \htext{$2$}
\move(-7.4 2.2) \htext{$\frac{1}{2}$}
\move(-5.75 1) \htext{$\frac{1}{2}$}
\move(-3.95 -0.3) \htext{$-b$}
\move(5.1 1) \htext{$\frac{\alpha_{1}}{\beta_{1}}$}
\move(6.2 2.2) \htext{$\frac{\alpha_{n}}{\beta_{n}}$}

\move(0 -4)

\end{texdraw}
\end{center}

\caption{Surgery presentation of $(\ns ; g \, | \, b;(\alpha_{1},\beta_{1}),\ldots,(\alpha_{n},\beta_{n}))$}\label{fig-A2}
\end{figure}

For completeness we will also state
the results in terms of the non-normalized
Seifert invariants due to W.\ D.\ Neumann, see \cite{JankinsNeumann}.
For a Seifert manifold $X$ with non-normalized Seifert invariants
$\{\ep;g;(\alpha_{1},\beta_{1}),\ldots,(\alpha_{n},\beta_{n})\}$
the invariants $\ep$ and $g$ are as above. The
$(\alpha_{j},\beta_{j})$ are here pairs of coprime integers with
$\alpha_{j} >0$ but not necessarily with $0< \beta_{j}<\alpha_{j}$.
These pairs are not invariants of $X$, but can be varied
according to certain rules. 
In fact, $X$ has a surgery presentation as shown in
Fig.~\ref{fig-A1} with $b=0$ if $\ep=\os$ and
as shown in Fig.~\ref{fig-A2} with $b=0$ if $\ep=\ns$.
The Seifert Euler number of $X$ is $-\sum_{j=1}^{n} \beta_{j}/\alpha_{j}$
(which is an invariant of the Seifert fibration $X$).
For more details, see 
\cite[Sect.~I.1]{JankinsNeumann}.

\section{Quantum invariants of Seifert manifolds}\label{sec-Quantum-invariants}

In this section we explain our calculation of 
the $\frg$--invariants of the Seifert manifolds, $\frg$ being
an arbitrary complex finite dimensional simple Lie algebra.
Our starting point is a formula
for the RT--invariants associated to an arbitrary modular category of
the Seifert manifolds, derived in \cite{Hansen1}. 

\rk{The RT--invariants of the Seifert manifolds for modular categories}
Let us first give some preliminary remarks on modular categories. We use
notation as in \cite{Turaev}.
Let $\left( \mV, \{ V_{i} \}_{i \in I } \right)$ 
be an arbitrary modular category with braiding $c$ and
twist $\theta$. The ground ring is $K=\Hom_{\mV}(\Io,\Io)$,
where $\Io$ is the unit object.
Let $i \mapsto i^{*}$ be the involution in $I$
determined by the condition that $V_{i^{*}}$ is isomorphic to
the dual of $V_{i}$. 
An element $i \in I$ is called self-dual if $i=i^{*}$. For such an element we
have a $K$--module isomorphism $\Hom_{\mV}(V \otimes V, \Io ) \cong K$, $V=V_{i}$.
The map $x \mapsto x(\id_{V} \otimes \theta_{V})c_{V,V}$ is a $K$--module
endomorphism of $\Hom_{\mV}(V \otimes V, \Io )$, so is a multiplication by a certain
$\vep_{i} \in K$. By the definition of the braiding and twist we have
$(\vep_{i})^{2} =1$. In particular $\vep_{i} \in \{ \pm 1 \}$ if $K$ is a field.
There is a distinguished element in $I$ denoted $0$, such that $V_{0}=\Io$.

The $S$-- and $T$--matrices of $\mV$ are the matrices $S=(S_{ij})_{i,j \in I}$,
$T=(T_{ij})_{i,j \in I}$ given by 
$S_{ij}=\tr(c_{V_{j},V_{i}} \circ c_{V_{i},V_{j}})$ and
$T_{ij}=\delta_{ij}v_{i}$, where $\tr$ is the categorical trace of $\mV$,
$\delta_{ij}$ is the Kronecher delta equal to $1$ if $i=j$ and zero elsewhere,
and $v_{i} \in K$ such that $\theta_{V_{i}}=v_{i}\id_{V_{i}}$.

Assume that $\mV$ has a rank $\mD$, i.e.\ an element of $K$ satisfying
$$
\mD^{2}=\sum_{i \in I} \dim(i)^{2},
$$
where $\dim(i)=\dim(V_{i})=\tr(\id_{V_{i}})$. We let
$$
\Delta = \sum_{i \in I} v_{i}^{-1} \dim(i)^{2}.
$$
Moreover, let $\tau=\tau_{(\mV,\mD)}$ be the RT--invariant
associated to  $\left( \mV, \{ V_{i} \}_{i \in I },\mD\right)$, cf.\ \cite[Sect.~II.2]{Turaev}.
For a tuple of integers $\mC =(m_{1},\ldots,m_{t})$, let
$$
G^{\mC} = T^{m_{t}}S\cdots T^{m_{1}}S.
$$
The Rademacher Phi function is 
defined on $PSL(2,\Z)=SL(2,\Z)/\{\pm 1\}$ by
\begin{equation}\label{eq:Phi}
\Phi \left( \left[ \begin{array}{cc}
                  p & r \\
                  q & s
                  \end{array}
\right] \right) = \left\{ \begin{array}{ll}
                  \frac{p+s}{q} - 12(\sign (q))\s (s,|q|) & ,q \neq 0, \\
                  \frac{r}{s} & ,q=0,
                  \end{array}
\right.
\end{equation}
see \cite{RademacherGrosswald}.
Here, for $q \neq 0$, the Dedekind sum $\s (s,q)$ is given by
\begin{equation}\label{eq:dedekindsum}
\s (s,q)= \frac{1}{4|q|} \sum_{j=1}^{|q|-1} \cot\frac{\pi j}{q} \cot \frac{\pi s j}{q}
\end{equation}
for $|q|>1$ and $\s (s,\pm 1)=0$.
We put $a_{\os}=2$ and $a_{\ns}=1$. Moreover, let
$b_{j}^{(\os)}=1$ and
$b_{j}^{(\ns)}=\delta_{j,j^{*}}$, $j \in I$.
Given pairs $(\alpha_{j},\beta_{j})$ of coprime integers we
let $\mC_{j}=(a_{1}^{(j)},\ldots,a_{m_{j}}^{(j)})$ be a continued
fraction expansion of $\alpha_{j}/\beta_{j}$,
$j=1,2,\ldots,n$, i.e.\
$$
\frac{\alpha_{j}}{\beta_{j}}=a_{m_{j}}^{(j)}-\frac{1}{a_{m_{j}-1}^{(j)}-\dfrac{1}{\cdots -\dfrac{1}{a_{1}^{(j)}}}}.
$$

\begin{thm}\label{invariants}{\rm\cite{Hansen1}}\qua
Let
$M=(\ep;g\;|\;b;(\alpha_{1},\beta_{1}),\ldots,(\alpha_{n},\beta_{n}))$, $\ep=\os,\ns$.
Then
\begin{eqnarray*}
&&\tau(M) = (\Delta\mD^{-1})^{\sigma_{\ep}} \mD^{a_{\ep}g-2-\sum_{j=1}^{n} m_{j}} \\
 && \hspace{1.0in} \times \sum_{j \in I} \left(\vep_{j}\right)^{a_{\ep}g} b_{j}^{(\ep)} v_{j}^{-b} \dim(j)^{2-n-a_{\ep}g} \left( \prod_{i=1}^{n} (SG^{\mC_{i}})_{j,0} \right),
\end{eqnarray*}
where
$$
\sigma_{\ep}=(a_{\ep}-1)\sign(E) + \sum_{j=1}^{n} \sign(\alpha_{j}\beta_{j}) + \frac{1}{3} \sum_{j=1}^{n} \left( \sum_{k=1}^{m_{j}} a_{k}^{(j)} -\Phi(B^{\mC_{j}}) \right).
$$
Here $E=-\left( b+\sum_{j=1}^{n} \frac{\beta_{j}}{\alpha_{j}} \right)$ is the Seifert Euler
number. 

The RT--invariant $\tau(M)$ of the Seifert manifold $M$ with
non-normalized Seifert invariants
$\{ \ep;g;(\alpha_{1},\beta_{1}),\ldots,(\alpha_{n},\beta_{n})\}$
is given by the same expression with the exceptions that the factor $v_{j}^{-b}$
has to be removed and $E=-\sum_{j=1}^{n} \frac{\beta_{j}}{\alpha_{j}}$.\HS
\end{thm}

The theorem is also valid in case $n=0$. In this case one just
has to put all sums $\sum_{j=1}^{n}$
equal to zero and all products $\prod_{i=1}^{n}$ equal to $1$.
Note that $\ep_{j}^{k}=1$ if $k$ is even and $\ep_{j}^{k}=\ep_{j}$
if $k$ is odd since $\ep_{j}^{2}=1$. In particular,
$\left(\vep_{j}\right)^{a_{\ep}g}=1$ if $\ep=\os$. The sum 
$\sum_{j=1}^{n} \sign(\alpha_{j}\beta_{j})$ is of course
equal to $n$ for normalized Seifert invariants. 

Let us next consider the lens spaces. For $p,q$ a pair of
coprime integers, recall that $L(p,q)$ is given by surgery
on $S^{3}$ along the unknot with surgery coefficient $-p/q$.
In the following corollary we include the possibilities
$L(0,1)=S^{1} \times S^{2}$ and $L(1,q)=S^{3}$, $q \in \Z$.

\begin{cor}\label{cor-lens-spaces}{\rm\cite{Hansen1}}\qua
Let $p,q$ be a pair of coprime integers. If $q \neq 0$ we let
$(a_{1},\ldots,a_{m-1})$ be
a continued fraction expansion of $-p/q$.
If $q=0$, put $m=3$ and $a_{1}=a_{2}=0$. 
Then
$$
\tau(L(p,q))=(\Delta\mD^{-1})^{\sigma} \mD^{-m} G^{\mC}_{0,0},
$$
where $\mC=(a_{1},\ldots,a_{m-1},0)$ and
$\sigma=\frac{1}{3} \left( \sum_{l=1}^{m-1} a_{l} - \Phi(B^{\mC}) \right)$.\HS 
\end{cor}

\rk{The RT--invariants of the Seifert manifolds for the classical Lie algebras}
It is a well-known fact that quantum deformations of the classical Lie
algebras at roots of unity induce modular categories, see \cite{Kirillov},
\cite{BakalovKirillov}, \cite{Le}. Let us provide the details needed.
For simplicity we will only consider simply laced Lie algebras in this paper,
except in \refrem{non-simply} where we give a few remarks
with respect to what have to be adjusted
to include the general case. (See also \cite{HansenTakata} for the general case.)
Therefore, let in the following 
$\frg$ be a fixed complex finite dimensional simple and simply laced Lie algebra.

First let us fix some notation for $\frg$. Let $\frh$ be a Cartan subalgebra of 
$\frg$, and let $\alpha_{1},\ldots,\alpha_{l}$ be a set of simple (basis)
roots in the dual space of $\frh$.
We denote by $\frhR^{*}$ the $\R$--vector space spanned 
by $\alpha_1,\ldots,\alpha_l$ and let
$\la \;,\; \ra$ be the inner product on $\frhR^{*}$ defined by
$\la \alpha_{i},\alpha_{j} \ra=a_{ij}$, $\left( a_{ij} \right)_{1 \leq i,j \leq l}$
being the Cartan matrix for $\frg$. In particular, all roots have length $\sqrt{2}$.
The root lattice $\rl$ is the $\Z$--lattice generated by 
$\alpha_1,\ldots,\alpha_l$, and the weight lattice $\wl$ is the
$\Z$--lattice generated by the fundamental weights
$\lambda_{1},\ldots,\lambda_{l}$,
i.e.\ $\lambda_{i} \in \frhR^{*}$ such that
$\la \lambda_{i},\alpha_{j} \ra=\delta_{ij}$ for all $i,j \in \{1,2,\ldots,l\}$. 
The (open) fundamental Weyl chamber is the set
$$
C=\{ x \in \frhR^{*} \; | \; \langle x,\alpha_{i} \rangle > 0, i=1,\ldots,l\}.
$$
For a positive integer $k$, the $k$--alcove is the (closed) set 
$$
C_{k}=\{ x \in \bar{C} \; | \; \langle x,\alpha_{0} \rangle \leq k\},
$$
where $\bar{C}$ is the topological closure of $C$ and
$\alpha_{0}$ is the highest root of $\frg$, i.e.\ $\alpha_{0}$ is the
unique root in $C$. The Weyl group is denoted $W$.

Let $q=e^{\pi \I/r}$, where $r$ is an integer $\geq h^{\vee}$.
Here $h^{\vee}$ is the dual Coxeter number of $\frg$ (equal to the Coxeter
number $h$ of $\frg$, since $\frg$ is simply laced).
By $U_{q}(\frg)$ we denote the quantum group associated to these data
as defined by Lusztig, see \cite[Part V]{Lusztig}. We follow
\cite[Sect.~1.3 and 3.3]{BakalovKirillov} here but will mostly
use notation from \cite{Turaev} for modular categories as above.
(Note that what we denote
$U_{q}(\frg)$ here is denoted $U_{q}(\frg)|_{q=e^{\pi\I/r}}$
in \cite{BakalovKirillov}.)
Let $\left(\mV_{r}^{\frg}, \{ V_{i} \}_{i \in I} \right)$ be the modular category
induced by the representation theory of $U_{q}(\frg)$, cf.\ \cite[Theorem 3.3.20]{BakalovKirillov}.
In particular, the index set for the simple objects
is $I=\inte(C_{r}) \cap \wl$.
We use here the shifted indexes (shifted by $\rho$) (contrary to \cite{BakalovKirillov}).
(Normally the irreducible modules
of $U_{q}(\frg)$ (of type $1$), $q$ a formal variable, are indexed by
the cone of dominant integer weights $\wl_{+}$. Here we
denote the irreducible module associated to $\mu \in \wl_{+}$
by $V_{\mu + \rho}$.) For $q$ a root of unity as above, $V_{\lambda}$
is an irreducible module of $U_{q}(\frg)$ of non-zero dimension if
$\lambda \in I$.
The involution $I \rightarrow I$, $\lambda \mapsto \lambda^{*}$,
is given by
$\lambda^{*} = -w_{0}(\lambda-\rho)+\rho$, where $w_{0}$
is the longest element in $W$ and $\rho$ is
half the sum of positive roots. The distinguished element $0 \in I$
is equal to $\rho$. According to \cite[Theorem 3.3.20]{BakalovKirillov}
we can use
\begin{equation}\label{eq:rank}
\mD=r^{l/2} \left| \frac{\vol (\rl) }{\vol (\wl) }\right|^{1/2}
  \left( \prod_{\alpha \in \Delta_{+}} 2\sin \left( \frac{\pi\la \alpha,\rho \ra }{r} \right) \right)^{-1}
\end{equation}
as a rank of $\mV_{r}^{\frg}$.
Here $\Delta_{+}$
is the set of positive roots.
According to the same theorem we have
\begin{equation}\label{eq:anomaly}
\Delta \mD^{-1} = \omega^{-3},
\end{equation}
where
\begin{equation}\label{eq:omega}
\omega = e^{\frac{2\pi\I c}{24}} = \exp\left( \frac{\pi\I}{h} |\rho|^{2} \right) \exp\left( -\frac{\pi\I}{r} |\rho|^{2} \right),
\end{equation}
where $c=\frac{r-h}{r}\dim(\frg)$ is the central
charge. The last equality in (\ref{eq:omega}) follows from Freudenthal's strange formula
$|\rho|^{2}/h=\dim \frg/12$.

The matrices $S$ and $T$ for $\mV_{r}^{\frg}$ are tightly related to a certain
unitary representation $\mR=\mR_{r}^{\frg}$ of $SL(2,\Z)$. On the standard
generators
\begin{equation}\label{eq:generators}
\Xi = \left( \begin{array}{cc}
		0 & -1 \\
		1 & 0
		\end{array}
\right),\hspace{.2in}
\Theta = \left( \begin{array}{cc}
                1 & 1 \\
		0 & 1
		\end{array}
\right) 
\end{equation}
of $SL(2,\Z)$ we have
\begin{eqnarray}\label{eq:mR}
&&\mR(\Xi)_{\lambda \mu} = \frac{\I^{\npr}}{r^{l/2}}
        \left| \frac{\vol (\wl) }{\vol (\rl) }\right|^{1/2} \sum_{w \in W} \det(w) 
        \exp \left( -\frac{2\pi\I}{r} \la w(\lambda), \mu\ra \right), \nonumber \\
&&\mR(\Theta)_{\lambda \mu} = \delta_{\lambda \mu}
        \exp \left( \frac{\pi\I}{r} \la \lambda,\lambda \ra
           -\frac{\pi\I}{h} \la \rho,\rho \ra \right)
\end{eqnarray}
for $\lambda, \mu \in I$.
In the following we also write $\tilde{U}$ for $\mR(U)$.
By using the results in \cite[Sect.~3.3]{BakalovKirillov},
in particular \cite[Theorem 3.3.20]{BakalovKirillov}, we find
\begin{equation}\label{eq:ST}
S_{\lambda \mu} = \mD \tilde{\Xi}_{\lambda \mu}, \hspace{.2in}T_{\lambda \mu} = \omega \tilde{\Theta}_{\lambda \mu}
\end{equation}
for $\lambda,\mu \in I$.
Let $\mC=(a_{1},\ldots,a_{n}) \in \Z^{n}$ and let $m \in \{0,1\}$. 
By (\ref{eq:ST}) we immediately get
\begin{equation}\label{eq:sumformula}
(S^{m}G^{\mC})_{\lambda \rho} = \mD^{m+n}\omega^{\sum_{j=1}^{n} a_{j}} \left( \tilde{\Xi}^{m} \tilde{\Theta}^{a_{n}}\tilde{\Xi}\tilde{\Theta}^{a_{n-1}}\cdots\tilde{\Theta}^{a_{1}}\tilde{\Xi}\right)_{\lambda \rho}
\end{equation}
for $\lambda \in I$. Finally we have for any $\lambda \in I$ that
$$
\dim(\lambda)=S_{\lambda \rho}=\mD \tilde{\Xi}_{\lambda \rho}=\mD r^{-l/2}
  \left| \frac{\vol (\wl) }{\vol (\rl) }\right|^{1/2}
  \prod_{\alpha \in \Delta_{+}} 2\sin \left( \frac{\pi\la \alpha,\lambda \ra }{r} \right),
$$
see also \cite[Formulas (3.3.2) and (3.3.5)]{BakalovKirillov}.
All the above data can now be put into the expression in \refthm{invariants}
to give a formula for $\tau_{r}^{\frg}(X)$,
$X$ an arbitrary Seifert manifold, where $\tau_{r}^{\frg}$ is the RT--invariant associated
to $\mV_{r}^{\frg}$. However, the formula to emerge is not detailed enough
to be of any use, at least not when it comes to a determination of asymptotics
of the invariants. The reason is, that the formula will contain matrix
products as in the right-hand side of (\ref{eq:sumformula}).

A way out of this problem is to determine nice formulas for the entries of
$\mR(U)$ in terms of the entries of $U$. This has in fact been done for
$\frg=\frsl_{2}(\C)$ by Jeffrey \cite{Jeffrey1}, \cite{Jeffrey2}.
Here method is to
write the matrix $U \in SL(2,\Z)$ as a product in the generators
$\Xi$ and $\Theta$ and make a certain induction argument. A main ingredient
is a reciprocity formula for Gaussian sums. We use
a similar argument to extend Jeffrey's results to an arbitrary
complex finite dimensional simple Lie algebra. The
following theorem generalizes
\cite[Propositions 2.7 and 2.8]{Jeffrey2}.

\begin{thm}\label{rep}
Let $U=\left( \begin{array}{cc}
		a & b \\
		c & d
		\end{array}
\right) \in SL(2,\Z)$ with $c \neq 0$. 
Then there exists an $\ep \in \{ \pm 1 \}$ such that 
\begin{eqnarray*}
&&\mR(\ep U)_{\lambda \mu} = \frac{ \I^{\npr} \sign(\ep c)^{\npr}}{(r|c|)^{l/2} \vol (\rl)}
      \exp \left( -\frac{\pi \I}{h} |\rho|^{2} \Phi(U) \right)\\
 && \hspace{1.0in} \times \exp \left( \frac{\pi \I}{r} \frac{d}{c} |\mu|^{2} \right) \sum_{\nu \in \rl/c \rl} 
      \exp \left( \frac{\pi \I} {r} \frac{a}{c} |\lambda+r\nu|^{2} \right) \\
 && \hspace{1.0in} \times \sum_{w\in W} \det (w)
      \exp\left( -\frac{2\pi\I}{r\ep c} \la \lambda + r\nu,w(\mu) \ra \right). 
\end{eqnarray*}
\end{thm}

The function $\Phi$ is given in (\ref{eq:Phi}).
If $c=0$, then $U= \ep \Theta^{b}$ for some $b \in \Z$ and $\ep \in \{ \pm 1 \}$ and
the expression for $\mR(\ep U)_{\lambda \mu}$ follows immediately from (\ref{eq:mR}).
At first sight the above theorem looks a little strange because of the
undetermined sign $\ep$. This sign has to do with the fact
that $\mR$ is a representation of $SL(2,\Z)$ and not
of $PSL(2,\Z)$ (except for $\frg=\frsl_{2}(\C)$, where $\mR$
is in fact a representation of $PSL(2,\Z)$).
However, as we shall see now, we will get rid of this sign
in the cases we need.
In fact, according to \refthm{invariants} and (\ref{eq:sumformula}), we only need the expression for
$\mR(U)_{\lambda \mu}$ in case $\lambda$
or $\mu$ is equal to $\rho$
in the calculation of the invariants of the Seifert manifolds.
Since $\rho^{*}=\rho$, $\Xi^{2}=-1$, and
$\mR(\Xi^{2})_{\lambda \mu}=\delta_{\lambda \mu^{*}}$,
we have
$$
\mR(-U)_{\lambda \rho} = \mR(U)_{\lambda \rho},\hspace{.2in} \mR(-U)_{\rho \lambda} = \mR(U)_{\rho \lambda}
$$
for all $\lambda \in I$. By using this fact
and the Weyl denominator formula one can show the following corollary
to \refthm{rep}.
(To show the first formula in \refcor{cor:rho} one also has to use unitarity of
$\mR$.)

\begin{cor}\label{cor:rho}
Let $U=\left( \begin{array}{cc}
		a & b \\
		c & d
		\end{array}
\right) \in SL(2,\Z)$ with $c \neq 0$. Then 
\begin{eqnarray*}
&&\mR(U)_{\lambda \rho} = \frac{ \I^{\npr} \sign(c)^{\npr}}{(r|c|)^{l/2} \vol (\rl)}
      \exp \left( -\frac{\pi \I}{h} |\rho|^{2} \Phi(U) \right) \\
 && \hspace{0.8in} \times \exp \left( \frac{\pi \I}{r} \frac{a}{c} |\lambda|^{2} \right) \sum_{\nu \in \rl/c \rl} 
      \exp \left( \frac{\pi \I}{r} \frac{d}{c} |\rho+r\nu|^{2} \right) \\
 && \hspace{0.8in} \times \sum_{w\in W} \det (w)
      \exp\left( -\frac{2\pi\I}{rc}
         \la \rho +r\nu,w(\lambda) \ra \right).
\end{eqnarray*}
If $a \neq 0$ we also have
\begin{eqnarray*}
&&\mR(U)_{\lambda \rho} = \frac{ \I^{\npr} \sign(c)^{\npr}}{(r|c|)^{l/2} \vol (\rl)} \exp \left( \frac{\pi \I}{r} \frac{b}{a} |\rho|^{2} \right) \\
 && \hspace{0.8in} \times\exp \left( -\frac{\pi \I}{h} |\rho|^{2} \Phi(U) \right) \\
 && \hspace{0.8in} \times \sum_{w\in W} \det (w) 
      \sum_{\nu \in \rl/c \rl}
      \exp\left( \frac{\pi\I}{r} \frac{a}{c} |\lambda +r\nu-\frac {w(\rho)}{a}|^{2} \right).  
\end{eqnarray*}\HS
\end{cor}

Because of the length and technical
nature of the proof of \refthm{rep}, we defer the argument to the Appendix,
and will only give the main ideas there. Detailed arguments will appear in \cite{HansenTakata}.

Given a pair of coprime integers $(\alpha,\beta)$,
$\alpha>0$, we let $\beta^{*}$ be the inverse of $\beta$ in the multiplicative
group of units in $\Z/\alpha\Z$. The following theorem is a generalization
of \cite[Theorem 8.4]{Hansen1} (which concerns the case $\frg=\frsl_{2}(\C)$).
The proof follows closely the proof of \cite[Teorem 8.4]{Hansen1} and is
therefore left out here. (One has to use the first formula in \refcor{cor:rho}.)

\begin{thm}\label{Lie-Seifert}
Let $M=(\ep;g\;|\;b;(\alpha_{1},\beta_{1}),\ldots,(\alpha_{n},\beta_{n}))$,
$\ep \in \{ \os, \ns\}$. Then
\begin{eqnarray*}
\tau_{r}^{\frg}(M) &=& \exp \left( \frac{\pi\I}{r}|\rho|^{2} \left[ 3(a_{\ep}-1)\sign(E) -E - 12\sum_{j=1}^{n} \s (\beta_{j},\alpha_{j}) \right] \right) \\
  &&\hspace{.1in} \times \frac{\I^{n\npr} r^{l(a_{\ep}g/2-1)}}{2^{\npr(n+a_{\ep}g-2)}\vol(\rl)^{2-a_{\ep}g}}
  \frac{1}{\mA^{l/2}} e^{\frac{3\pi\I}{h}|\rho|^{2}(1-a_{\ep})\sign(E)} Z_{\ep}^{\frg}(M;r),
\end{eqnarray*}
where $\s (\beta_{j},\alpha_{j})$ is given by {\em (\ref{eq:dedekindsum})}, $\mA=\prod_{j=1}^{n} \alpha_{j}$, and
\begin{eqnarray*}
Z_{\ep}^{\frg}(M;r) &=& \sum_{\lambda \in I} b_{\lambda}^{(\ep)} \vep_{\lambda}^{a_{\ep}g}
\left( \prod_{\alpha \in \Delta_{+}} 
\sin ^{2-n-a_{\ep}g} \left( \frac{\pi\la \lambda,\alpha \ra}{r} \right)\right)
\exp \left( \frac{\pi\I}{r} E |\lambda|^{2} \right) \\
 & & \hspace{.2in} \times \sum_{w_{1},\ldots,w_{n}\in W} 
  \sum_{\nu_{1} \in \rl/\alpha_{1} \rl} \ldots \sum_{\nu_{n} \in \rl/\alpha_{n} \rl}
  \left(\prod_{j=1}^{n} \det(w_{j})\right) \\
 & & \hspace{0.6in} \times
     \exp \left( -\pi \I \sum_{j=1}^{n} \frac{\beta_{j}^{*}}{\alpha_{j}} \left( r|\nu_{j}|^{2} 
      + 2\la w_{j}(\rho),\nu_{j} \ra \right) \right) \\
 & & \hspace{0.6in} \times  
     \exp \left( - \frac{2 \pi\I}{r} \la \lambda,  \sum_{j=1}^{n} \frac{ r\nu_{j} + w_{j}(\rho)}{\alpha_{j}} \ra \right).
\end{eqnarray*}
The RT--invariant $\tau_{r}^{\frg}(M)$ of the Seifert manifold $M$
with non-normalized Seifert invariants
$\{\ep;g;(\alpha_{1},\beta_{1}),\ldots,$ $(\alpha_{n},\beta_{n})\}$
is given by the same expression.\HS
\end{thm}

The theorem is also valid in case $n=0$. In this case one
just has to put the sum
$\sum_{w_{1},\ldots,w_{n}\in W} 
\sum_{\nu_{1} \in \rl/\alpha_{1} \rl} \ldots \sum_{\nu_{n} \in \rl/\alpha_{n} \rl}$
in $Z_{\ep}(M;r)$ equal to $1$,$\ep=\os,\ns$,
and put $\mA=1$ and 
$\sum_{j=1}^{n}\s (\beta_{j},\alpha_{j})=0$. 

Let us finally consider the lens space $L(p,q)$.
Let $b,d$ be any integers such that
$U=\left( \begin{array}{cc}
                q & b \\
		p & d
		\end{array}
\right) \in SL(2,\Z)$. Assume $q \neq 0$, let $V=-\Xi U =\left( \begin{array}{cc}
                p & d \\
		-q & -b
		\end{array}
\right)$, and let $C'=(a_{1},a_{2},\ldots,a_{m-1}) \in \Z^{m-1}$ such that
$B^{\mC'}=V$. Then $\mC'$ is a continued fraction expansion of $-p/q$ and
$U=\Xi V=B^{\mC}$, where $\mC=(a_{1},a_{2},\ldots,a_{m-1},0)$.
By \refcor{cor-lens-spaces}, (\ref{eq:anomaly}) and (\ref{eq:sumformula})
we therefore get
\begin{equation}\label{eq:RTlens}
\tau_{r}^{\frg}(L(p,q))= \omega ^{\Phi(U)} \tilde{U}_{\rho\rho},
\end{equation}
where $\omega$ is given by (\ref{eq:omega}).
If $q=0$ we have $p=1$ and $L(p,q)=S^{3}$. In this case we have
$\tau_{r}^{\frg}(L(p,q))=\mD^{-1}$. We also have
$U=\Xi \Theta^{d}$, so by using
(\ref{eq:mR}), (\ref{eq:rank}) and (\ref{eq:omega})
we find that the right-hand side of (\ref{eq:RTlens})
is also equal to $\mD^{-1}$.
The identity (\ref{eq:RTlens}) coincides with \cite[Formula (3.7)]{Jeffrey2}
for $\frg=\frsl_{2}(\C)$, see also \cite[Formula (49)]{Hansen1}.

\begin{rem}\label{non-simply}
Let us briefly mention the adjustments to be done for including the non-simply
laced Lie algebras. In the general case the root of unity $q=e^{\pi \I/(dr)}$,
where $d=1$ for $\frg$ simply laced, $d=2$ if $\frg$ belongs to the series
$BCF$ and $d=3$ if $\frg$ is of type $G_{2}$. Moreover, $\rl$ is in general
the coroot lattice, which is dual to the weight lattice.
The inner product $\la \; , \; \ra$ in $\frhR^{*}$ is induced by an invariant
bilinear form on $\frg$, and is normalized such that
a long root has length $\sqrt{2}$. We stress that $\alpha_{0}$ is
the long highest root of $\frg$.
\end{rem}

\section{The asymptotic expansion conjecture and Seifert manifolds}\label{sec-The-asymptotic}

\noindent For $X$ a fixed closed oriented $3$--manifold 
we consider $r \mapsto \tau_{r}^{\frg}(X)$
as a complex valued function on $\{h^{\vee},h^{\vee}+1,h^{\vee}+2,\ldots\}$. We
are interested in the behaviour of this function in the limit of large $r$, i.e.
$r \rightarrow \infty$.

It is believed that Witten's TQFT's associated to $G$
coincides with the TQFT's of Reshetikhin--Turaev associated to $\frg$, where
$G$ is a simply connected compact simple Lie group with complexified Lie
algebra $\frg$.
In particular it is conjectured that Witten's leading large $k$ asymptotics for 
$Z_{k}^{G}(X)$
should be valid for the function $r \mapsto \tau_{r}^{\frg}(X)$ in the 
limit $r \ria \infty$, and furthermore, that this function should have a full 
asymptotic expansion.
The precise formulation of this is stated in the following conjecture,
called the asymptotic expansion conjecture (AEC).

\begin{conj}[J.\ E.\ Andersen \cite{Andersen1}, \cite{Andersen2}]\label{AEC} 
Let $\{\alpha_{1},\ldots,\alpha_{M} \}$
be the set
of values of the Chern--Simons functional of flat $G$ connections on 
a closed oriented $3$--manifold $X$. Then there exist $d_{j} \in \Q$,
$\tilde{I}_{j} \in \Q / \Z$, $b_{j} \in \R_{+}$ and $c_{m}^{j} \in \C$
for $j=1,\ldots,M$ and $m=1,2,3,\ldots$ such that
\begin{equation}\label{eq:AEC}
\tau_{r}^{\frg}(X) \sim_{r \ria \infty} \sum_{j=1}^{M} b_{j}e^{2\pi \sqrt{-1} r \alpha_{j}} r^{d_{j}} e^{\pi \sqrt{-1} \tilde{I}_{j}/4} \left(1+\sum_{m=1}^{\infty} c_{m}^{j} r^{-m} \right),
\end{equation}
that is, for all $N=0,1,2,\ldots$
$$
\tau_{r}^{\frg}(X) = \sum_{j=1}^{M} b_{j}e^{2\pi \sqrt{-1} r \alpha_{j}} r^{d_{j}} e^{\pi \sqrt{-1} \tilde{I}_{j}/4} \left(1+\sum_{m=1}^{N} c_{m}^{j} r^{-m} \right) + o(r^{d-N})
$$
in the limit $r \ria \infty$, where $d=\max\{d_{1},\ldots,d_{M} \}$.
\end{conj}

As noticed by Andersen \cite{Andersen1}, \cite{Andersen2}, a
complex function defined on the positive integers has
at most one asymptotic expansion on the form (\ref{eq:AEC})
if the $\alpha_{j}$'s are rational (and mutually different),
see also \cite{Hansen2}.
This means that if the AEC is true and if we futhermore
have that the Chern--Simons invariants are rational 
(as conjectured by e.g.\ Auckly \cite{Auckly}), then all the quantities
$c_{m}^{j}$, $\alpha_{j}$, $d_{j}$,... are
topological invariants. 
The AEC was first proved for the mapping tori of finite order
diffeomorphisms of orientable surfaces of genus at least $2$ and for any $\frg$
by Andersen \cite{Andersen1} using the gauge theory definition of the invariants.
Note that these mapping tori are Seifert manifolds with orientable base
and Seifert Euler number equal to zero.
Later on, the AEC was proved for all Seifert manifolds with
orientable base or non-orientable base with even genus
in case $\frg=\frsl_{2}(\C)$, cf.\ \cite{Hansen2}.
The proof of this result is partly based on calculations of
Rozansky \cite{Rozansky}.
For more details about the AEC and conjectures about the topological
interpretation of the different parts of the asymptotic formula (\ref{eq:AEC})
we refer to \cite{Andersen1}, \cite{Andersen2},
\cite{Hansen2}. One can also find a review about the status of the AEC in
these references.

By elaborating on the expression (\ref{eq:RTlens})
along the same lines as in \cite[Sect.~3]{Jeffrey2}
(using the last formula in \refcor{cor:rho}) we find
the following generalization of \cite[Theorem 3.4]{Jeffrey2}
(valid in case $p \neq 0$).

\begin{thm}\label{asymp-lens}
The RT--invariant associated to $\frg$ of the lens space $L(p,q)$
is given by
\begin{gather}
 \tau_{r}^{\frg}(L(p,q)) = \frac{\sign(p)^{\npr}\I^{\npr}}{(r|p|)^{l/2} \vol(\rl)}
       \exp \left( \frac{\pi \I}{r} 12 \sign(p)\s (q,|p|) |\rho|^{2} \right)\notag \\
  \hspace{1.0in} \times \sum_{w \in W} \det(w) 
       \exp \left( -\frac{2\pi \I}{pr} \la \rho, w(\rho) \ra \right)\notag \\
  \hspace{0.6in} \times \sum_{\nu \in \rl/p\rl}
       \exp \left( \pi \I \frac{q}{p} r |\nu|^{2} \right)
       \exp \left( 2\pi \I \frac{1}{p} \la \nu, q\rho - w(\rho) \ra \right).\tag*{\qed}
\end{gather}
\end{thm}

From this theorem
it is obvious, that the large $r$ asymptotics
of $\tau_{r}^{\frg}(L(p,q))$ is on the same form as in (\ref{eq:AEC})
(expand the factor
$\exp \left( \frac{\pi \I}{r} 12 \sign(p)\s (q,|p|) |\rho|^{2} \right)$ $\exp \left( -\frac{2\pi \I}{pr} \la \rho, w(\rho) \ra \right)$
as a power series in $r^{-1}$).
Proving the following conjecture will therefore finalize the proof of the AEC
for the invariants $\tau_{r}^{\frg}(L(p,q))$.

\begin{conj}
The set of values of the Chern--Simons functional of flat $G$ connections on
$L(p,q)$ is given by
$$
\left\{ \left. \; \frac{q}{2p} |\nu|^{2} \pmod{\Z} \;\; \right| \;\; \nu \in \rl/p\rl \; \right\}.
$$ 
\end{conj}

For $\frg = \frsl_{n}(\C)$ ($G=SU(n)$) this conjecture should follow from results in
\cite{Nishi}.

\section{Appendix. The proof of \refthm{rep}}\label{sec-Appendix}

In this section we will explain the ideas behind the
proof of \refthm{rep}. The underlying Lie algebra $\frg$ is still assumed
to be simply laced for simplicity. However, with some minor adjustments
the arguments given are also
true for the non-simply laced Lie algebras, see
\refrem{non-simply}.

The proof of \refthm{rep} builds mainly on the
key-lemma, \reflem{lem:main}. Let us introduce some notation.
For a tuple of integers $\mC=(m_{1},\ldots,m_{t})$, let
\begin{equation}\label{eq:Bmatrix}
B_{k}^{\mC} = \left( \begin{array}{cc}
                a_{k}^{\mC} & b_{k}^{\mC} \\
		c_{k}^{\mC} & d_{k}^{\mC}
		\end{array}
\right) = \Theta^{m_{k}}\Xi\Theta^{m_{k-1}}\Xi \ldots \Theta^{m_{1}}\Xi 
\end{equation}
for $k=1,2,\ldots,t$,
and let $B^{\mC}=B_{t}^{\mC}$, where $\Xi$ and $\Theta$ are given by
(\ref{eq:generators}). Moreover, we put
$$ 
a_{0}^{\mC}=d_{0}^{\mC}=1, \qquad b_{0}^{\mC}=c_{0}^{\mC}=0.
$$
We say that $\mC$ has length $|\mC|=t$.
If it is clear from the context what $\mC$ is
we write $a_{k}$ for $a_{k}^{\mC}$ etc.
From \cite[Proposition 2.5]{Jeffrey2}, the elements $a_{i},b_{i},c_{i},d_{i}$
satisfy the recurrence relations
\begin{align}\label{eq:abcd}
a_{k}=m_{k} a_{k-1}-c_{k-1},& \quad c_{k}=a_{k-1}, \\
b_{k}=m_{k} b_{k-1}-d_{k-1},& \quad d_{k}=b_{k-1} \nonumber
\end{align}
for $k=1,2,\ldots,t$.
One should note that the expressions (\ref{eq:mR}) for the entries of $\mR(\Xi)$
and $\mR(\Theta)$ are well-defined for all $\lambda, \mu \in \wl$.
Note also that if $\lambda$ or $\mu$ are elements of $\wl$
belonging to the boundary of $C_{r}$ then $\mR(\Xi)_{\lambda \mu}=0$.
This observation allows us to shift between $I=\inte(C_{r}) \cap \wl$
and $C_{r} \cap \wl$ as summation index set in formulas below.
This shift is important in the proof of \reflem{lem:main}.
Following Jeffrey \cite[Sect.~2]{Jeffrey1}, \cite[Sect.~2]{Jeffrey2} we consider
$$
\mT^{\mC}_{\lambda_{0},\lambda_{t+1}} = 
\sum_{\lambda_{1},\ldots,\lambda_{t} \in C_{r} \cap \wl}
\tilde{\Xi}_{\lambda_{t+1} \lambda_{t}} 
\tilde{\Theta}_{\lambda_t}^{m_t} \tilde{\Xi}_{\lambda_{t} \lambda_{t-1}} 
\tilde{\Theta}_{\lambda_{t-1}}^{m_{t-1}} \tilde{\Xi}_{\lambda_{t-1} \lambda_{t-2}} 
\cdots
\tilde{\Theta}_{\lambda_{1}}^{m_{1}} \tilde{\Xi}_{\lambda_{1} \lambda_{0}}
$$
for $\lambda_{0},\lambda_{t+1} \in C_{r} \cap \wl$, where we write
$\tilde{\Theta}_{\lambda}$ for $\tilde{\Theta}_{\lambda \lambda}$.
Then we have the following generalization of \cite[Lemma 2.6]{Jeffrey2}:

\begin{lem}\label{lem:main}
Assume that $\mC=(m_{1},\ldots,m_{t})$ is a sequence of integers such that
$a_{k}$ is nonzero for $k=1,\ldots,t$. Then
$$
\mT^{\mC}_{\lambda_{0},\lambda_{t+1}} = K^{\mC}_{\lambda_{0}} \sum_{w\in W} \det(w) 
      \sum_{\mu \in \rl/a_{t} \rl}
       \exp\left( -\frac {\pi\I c_{t}}{a_{t} r}
       \left|\lambda_{t+1}+r\mu+\frac {w(\lambda_{0})}{c_{t}}\right|^{2} \right),
$$
where 
\begin{eqnarray*}
K^{\mC}_{\lambda_{0}} &=& \frac{{\I}^{(t+1) \npr}}{(r|a_{t}|)^{l/2}\vol (\rl)}\; \zeta^{l \; D_{t}}
 \exp \left( -\frac{\pi\I}{h} (\sum_{i=1}^{t} m_{i}) |\rho |^{2} \right) \\ 
 & & \hspace{.3in} \times \exp \left( -\frac{\pi\I}{r} 
        (\sum_{i=1}^{t-1} \frac{1}{a_{i-1}a_{i}}) |\lambda_{0} |^{2} \right).
\end{eqnarray*}
Here $\zeta=\exp \frac{\pi \I}{4}$
and $D_{t}=\sign(a_{0}a_{1})+\cdots +\sign(a_{t-1}a_{t})$.\qed
\end{lem}

We will not give the proof of this lemma here, since it is long
and technical.
The lemma is proved by induction on the length of $\mC$. The
reciprocity formula for Gaussian sums, \refprop{prop:gauss}, plays a
prominent role in the proof. A proof of this reciprocity formula can be
found in \cite[Sect.~2]{Jeffrey1}.

Let $V$ be a real vector space of dimension $l$ with inner product 
$\langle \;, \;\rangle$, $\Lambda$ a lattice in $V$ and $\Lambda^{*}$ 
the dual lattice. 
For an integer $r$, a self-adjoint automorphism 
$B \co V \to V$, and an element  $\psi \in V$, 
we assume 
\begin{eqnarray*}
&&\frac{1}{2} \la \lambda,B r \lambda \ra, \;\;
\la \lambda,B \eta \ra, \;\;
r \la \lambda, \psi \ra \in \Z, \hspace{.2in}\forall \lambda, \eta \in \Lambda,\\
&&\frac{1}{2} \la \mu,B r \mu \ra, \;\;
\la \mu,r \xi \ra, \;\;
r \la \mu, \psi \ra \in \Z, \hspace{.2in}\forall \mu, \xi \in \Lambda^*\\
\end{eqnarray*}
and $B\Lambda^{*} \subseteq \Lambda^{*}$. Then we have the following:
 
\begin{prop}[Reciprocity formula for Gauss sums]\label{prop:gauss} 
\begin{gather}
\kern -1.4in \vol (\Lambda^{*}) \sum_{\lambda \in \Lambda / r\Lambda}
     \exp \left( \frac{\pi\I}{r} \langle \lambda,B \lambda\rangle\right) 
     \exp \left(2\pi\I \la \lambda,\psi \ra \right)\notag \\
\hspace{.2in} =\left( \det \frac B \I \right)^{-1/2} r^{l/2} 
     \sum_{\mu \in \Lambda^{*} / B\Lambda^{*}}
      \exp \left( -\pi r\I \la \mu+\psi, B^{-1}(\mu+\psi) \ra \right).\tag*{\qed}
\end{gather}
\end{prop}

In the proof of \reflem{lem:main} we use Proposition \ref{prop:gauss}
with $\Lambda =\wl$, the dual lattice being the
root lattice $\rl$. Basically we use the reciprocity formula
($t$ times, recursively)
to change the sums in the expression for $\mT^{\mC}$
to a sum with a range which does not depend on
$r$. Another main ingredient in the proof of \reflem{lem:main}
is symmetry considerations along the same lines as the discussion in
\cite[pp.~584--586]{Jeffrey2}, see in particular
\cite[Proposition 4.4]{Jeffrey2}.
In the proof of \reflem{lem:main} and \refthm{rep} we
use several results in
\cite[Sect.~2]{Jeffrey2}, in particularly \cite[Proposition 2.5]{Jeffrey2}.

\reflem{lem:main} nearly proves \refthm{rep}. There is, however, a small
hurdle to overcome because of the assumption on the $a_{k}$'s in the lemma.
The following small result does the job.

\begin{lem}\label{lem:matrixdecomposition}
Let $U=\begin{bmatrix}
a & b \\
c & d
\end{bmatrix} \in PSL(2,\Z)$ with $c \neq 0$. Then we can write
$U=V\Theta^{n}$,
where $n \in \Z$ and where $V$ is given in the following way:
If $a=0$, then $V =\Xi$; if $a\neq 0$, then
there exists a sequence of integers $\mC$ such that
$V=B^{\mC}$, see {\em (\ref{eq:Bmatrix})}, and such that $a^{\mC}_{k} \neq 0$,
$k=1,2,\ldots,|\mC|$.\HS
\end{lem}

From this lemma we see the origin of the undetermined sign $\ep$.
Let us use the above lemmas to sketch the proof of \refthm{rep}.

\begin{proof}[Proof of \refthm{rep}]
According to the previous lemma there exists an integer $n$, a sign $\ep \in \{ \pm 1 \}$,
and a $V \in SL(2,\Z)$ as in the \reflem{lem:matrixdecomposition}
such that $U=\ep V\Theta^{n}$.
First assume $a \neq 0$ and that $n=0$, i.e.\ assume that $\ep U=B^{\mC}$, where $\mC=(m_{1},\ldots,m_{t})$
and $a_{k} \neq 0$, $k=1,2,\ldots,t$. Let $\mC'=(m_{1},\ldots,m_{t-1})$.
Then, by \reflem{lem:main}, (\ref{eq:abcd}), and a small calculation, we get
\begin{eqnarray*}
&&\mR(\ep U)_{\lambda \mu} = \tilde{\Theta}_{\lambda \lambda}^{m_t} \mT^{\mC'}_{\mu,\lambda} \\
 && \hspace{.5in} = \sum_{w\in W} \det (w) 
      \sum_{\nu \in \rl/c \rl}
      \exp\left( \frac{\pi\I}{r} \frac{a}{c}
         \left|\lambda+r\nu-\frac {w(\mu)}{\ep a}\right|^{2} \right),        
\end{eqnarray*}
where
$$
K = K^{\mC'}_{\mu} \exp \left( -\frac{\pi \I}{h} m_{t} |\rho|^{2} \right) \exp \left( -\frac{\pi \I}{a_{t-1}a_{t-2}r} |\mu|^{2} \right) \exp \left( -\frac{\pi \I}{a_{t}a_{t-1}r} |\mu|^{2} \right).
$$ 
The formula for the entries of $\mR(\ep U)$
now follows by the fact $b/a+1/(ac)=d/c$
(except for the factor $K$ which we will not rewrite here).
Next assume that $\ep U = B^{\mC}\Theta^{n}$ with $n \neq 0$,
where $\mC$ is as above. Then
$$
B^{\mC} = \ep U \Theta^{-n} = \ep \left( \begin{array}{cc}
		a & b \\
		c & d
		\end{array}
\right)\left( \begin{array}{cc}
		1 & -n \\
		0 & 1
		\end{array}
\right)=\ep \left( \begin{array}{cc}
		a & -na+b \\
		c & -nc+d
		\end{array}
\right),
$$
and since \refthm{rep} is valid for $U=B^{\mC}$ (with $\ep =1$) we get
the result after a small calculation.
Finally one has to consider the case where $a=0$, in which case 
$\ep U=\Xi \Theta^{n}=
\left( \begin{array}{cc}
		0 & -1 \\
		1 & n
		\end{array}
\right)$. Here the result follows by inserting the expressions
in (\ref{eq:mR}) into
$\mR(\ep U)_{\lambda \mu}=\mR(\Xi)_{\lambda \mu}\mR(\Theta)^{n}_{\mu \mu}$.
\end{proof}

\Addresses\recd

\end{document}